\documentclass[11pt]{amsart}

\usepackage[usenames]{color}

\usepackage[colorlinks, linkcolor=reference,
            citecolor=citation,
          urlcolor=e-mail]{hyperref}

\definecolor{e-mail}{rgb}{0,.40,.80}
\definecolor{reference}{rgb}{.20,.60,.22}
\definecolor{citation}{rgb}{0,.40,.80}
\usepackage{amsmath,times,pslatex,enumitem}
\usepackage{amsthm}
\usepackage{amsfonts}
\usepackage{amssymb}
\usepackage{latexsym}

\usepackage{hyperref}

\usepackage[top=1.2in, bottom=1in, left=1.7in, right=1.2in]{geometry}

\newtheorem{example}{Example}

\newtheorem{definition}{Definition}

\newtheorem{remark}{Remark}
\newtheorem*{remark*}{Remark}

\newcommand{\Var}{\mathop{\mathrm{Var}}}

\newcommand{\Z}{\mathbb{Z}}

\def\P{{\mathbf{P}}}

\def\NP{{\mathbf{NP}}}

\def\SSP{{\mathbf{SSP}}}

\def\SMP{{\mathbf{SMP}}}

\def\KP{{\mathbf{KP}}}

\def\PCP{{\mathbf{PCP}}}

\def\GPCP{{\mathbf{GPCP}}}

\def\CA{{\mathcal{A}}}

\begin{document}

\title{Problems in group theory motivated by cryptography}

\author[]{Vladimir Shpilrain}
\address{Department of Mathematics, The City College of New York, New York,
NY 10031} \email{shpil@groups.sci.ccny.cuny.edu}
\thanks{Research of the author was partially supported by
the ONR (Office of Naval Research) grant N000141512164}

\begin{abstract}
This is a survey of algorithmic problems in group theory, old and
new, motivated by applications to cryptography.

\end{abstract}

\maketitle

\section{Introduction}

The object of this survey is to showcase algorithmic problems in
group theory motivated by (public key) cryptography.

In the core of most public key cryptographic primitives there is an
alleged practical irreversibility of some process, usually referred
to as a {\it one-way function with trapdoor}, which is a function
that is easy to compute in one direction, yet believed to be
difficult to compute the  inverse function on ``most" inputs without
special information, called the ``trapdoor". For example, the RSA
cryptosystem uses the fact that, while it is not hard to compute the
product of two large primes, to {\it factor} a very large integer
into its prime factors appears to be computationally hard. Another,
perhaps even more intuitively obvious,  example is that of the
function $f(x)=x^2$. It is rather easy to compute in many reasonable
(semi)groups, but the inverse function $\sqrt{x}$ is much less
friendly. This fact is exploited in Rabin's cryptosystem, with the
multiplicative semigroup of ${\Z_n}$ ($n$ composite) as the
platform. In both cases though, it is not immediately clear what the
trapdoor is. This is typically the most nontrivial part of a
cryptographic scheme.

For a rigorous definition of a one-way function we refer the reader
to \cite{Talbot}; here we just say that there should be an efficient
(which usually means polynomial-time with respect to the complexity
of an input) way to compute this function, but no visible
(probabilistic) polynomial-time algorithm for computing the inverse
function on ``most" inputs.

Before we get to the main subject of this survey, namely problems in
combinatorial and computational group theory motivated by
cryptography, we recall historically the first public-key
cryptographic scheme, the Diffie-Hellman key exchange protocol, to
put things in perspective. This is done in Section \ref{DH}. We note
that the platform group for the original Diffie-Hellman protocol was
finite cyclic. In Section \ref{ElGamal}, we show how to convert the
Diffie-Hellman  key exchange protocol to an encryption scheme, known
as the ElGamal cryptosystem.

In the subsequent sections, we showcase various problems about
infinite non-abelian groups. Complexity of these problems in
particular groups has been used in various cryptographic primitives
proposed over the last 20 years or so. We mention up front that a
significant shift in paradigm motivated by  research in cryptography
was moving to {\it search} versions of {\it decision} problems that
had been traditionally considered in combinatorial group theory, see
e.g. \cite{MSUbook2, witness}. In some cases,  decision problems
were used in cryptographic primitives (see e.g. \cite{Osin}) but
these occasions are quite rare.

The idea of using the complexity of infinite non-abelian groups in
cryptography goes back to Wagner and Magyarik  \cite{MW} who in 1985
devised a public-key protocol based on   the unsolvability of the
word problem for finitely presented groups (or so they thought).
Their protocol now looks somewhat naive, but it was pioneering. More
recently, there has been an increased interest in  applications of
non-abelian group theory to cryptography initially prompted by the
papers \cite{AAG, KLCHKP, sidel}.

We note that a separate question of interest that is outside of the
scope of this survey is what groups can be used as platforms for
cryptographic protocols. We refer the reader to the monographs
\cite{MSUbook1}, \cite{MSUbook2}, \cite{VMR_book} for relevant
discussions and examples; here we just mention that finding a
suitable platform (semi)group for one or another cryptographic
primitive is a challenging problem. This is currently an active area
of research; here we can mention that groups that have been
considered in this context include braid groups (more generally,
Artin groups), Thompson's group, Grigorchuk's group, small
cancellation groups, polycyclic groups, (free) metabelian groups,
various groups of matrices, semidirect products, etc.

Here is the list of algorithmic problems that we discuss in this
survey. In most cases, we consider search versions of the problems
as more relevant to cryptography, but there are notable exceptions.

\medskip

\noindent -- The word (decision) problem: Section \ref{WP}
\smallskip

\noindent -- The conjugacy problem: Section \ref{conjugacy}
\smallskip

\noindent -- The twisted conjugacy problem: Section \ref{twist}

\smallskip

\noindent -- The decomposition problem: Section \ref{decomposition}

\smallskip

\noindent -- The subgroup intersection problem: Section \ref{central}

\smallskip

\noindent -- The factorization problem: Section \ref{factorization}

\smallskip

\noindent -- The isomorphism inversion problem: Section \ref{isomorphism}

%

\smallskip

\noindent -- The subset sum and the knapsack problems: Section \ref{knapsack}

\smallskip

\noindent -- The Post correspondence problem: Section \ref{Post}

\smallskip

\noindent -- The hidden subgroup problem: Section \ref{HSP}

\medskip

\noindent Also, in Section \ref{semidirect} we show that using
semidirect products of (semi)groups as platforms for a
Diffie-Hellman-like key exchange protocol yields various peculiar
computational assumptions and, accordingly, peculiar search
problems.

In the concluding Section \ref{relations}, we describe relations
between some of the problems discussed in this survey.

\section{The Diffie-Hellman key exchange protocol}
\label{DH}

The whole area of public-key cryptography started  with the seminal
paper by Diffie and Hellman \cite{DH}. We quote from Wikipedia:
``Diffie-Hellman key agreement was invented in 1976 \dots and was
the first practical method for establishing a shared secret over an
unprotected communications channel." In 2002 \cite{Hsurvey}, Martin
Hellman gave credit to Merkle as well: ``The system \dots has since
become known as Diffie-Hellman key exchange. While that system was
first described in a paper by Diffie and me, it is a public-key
distribution system, a concept developed by Merkle, and hence should
be called `Diffie-Hellman-Merkle key exchange' if names are to be
associated with it. I hope this small pulpit might help in that
endeavor to recognize Merkle's equal contribution to the invention
of public-key cryptography."

U. S. Patent 4,200,770, now expired, describes the algorithm, and
credits Diffie, Hellman,  and Merkle as inventors.

The simplest, and original, implementation of the protocol uses the
multiplicative group $\Z_p^*$ of integers modulo $p$, where $p$ is
prime and $g$ is primitive mod $p$. A more general description of
the protocol uses an arbitrary finite cyclic group.

\begin{enumerate}

\item Alice and Bob agree on a finite cyclic group $G$ and a generating element $g$ in $G$.
 We will write the group $G$ multiplicatively.

\item Alice picks a random natural number $a$ and sends $g^a$ to Bob.

\item    Bob picks a random natural number $b$ and sends $g^b$ to Alice.

\item   Alice computes $K_A=(g^b)^a=g^{ba}$.

\item  Bob computes $K_B=(g^a)^b=g^{ab}$.
\end{enumerate}

Since $ab=ba$ (because $\Z$ is commutative), both Alice and Bob are
now in possession of the same group element $K=K_A= K_B$ which can
serve as the shared secret key.

The protocol is considered secure against eavesdroppers if $G$ and
$g$ are chosen properly. The eavesdropper, Eve, must solve the
 {\it Diffie-Hellman problem}
(recover $g^{ab}$ from $g^a$ and $g^b$) to obtain the shared secret
key. This is currently considered difficult for a ``good" choice of
parameters (see e.g. \cite{Menezes} for details).

An efficient algorithm to solve the  {\it discrete logarithm problem} (i.e., recovering $a$ from
$g$ and  $g^a$) would obviously solve the Diffie-Hellman problem,
making this and many other public-key cryptosystems insecure.
However, it is not known whether or not the discrete logarithm
problem is {\it equivalent} to the Diffie-Hellman problem.

We note that there is a ``brute force" method for solving the
discrete logarithm problem: the eavesdropper Eve can just go over
natural numbers $n$ from 1 up one at a time, compute $g^n$ and see
whether she has a match with the transmitted element. This will
require $O(|g|)$ multiplications, where  $|g|$ is the order of $g$.
Since in practical implementations $|g|$ is typically at least
$10^{300}$, this method is considered computationally infeasible.

This raises a question of computational efficiency for legitimate
parties: on the surface, it looks like legitimate parties, too, have
to perform $O(|g|)$ multiplications to compute $g^a$ or  $g^b$.
However, there is a faster way to compute $g^a$ for a particular $a$
by using the ``square-and-multiply" algorithm, based on the binary
form of $a$. For example, $g^{22}=(((g^2)^2)^2)^2 \cdot (g^2)^2
\cdot g^2$. Thus, to compute $g^a$, one actually needs $O(\log_2 a)$
multiplications, which is  feasible given the magnitude of $a$.

\subsection{The ElGamal cryptosystem} \label{ElGamal}

The  ElGamal cryptosystem \cite{ElGamal} is a public-key
cryptosystem which is based on the Diffie-Hellman key exchange. The
ElGamal protocol is used in the free GNU Privacy Guard software,
recent versions of PGP, and other cryptosystems. The Digital
Signature Algorithm (DSA) is a variant of the ElGamal signature
scheme, which should not be confused with the ElGamal encryption
protocol that we describe below.

\begin{enumerate}

\item Alice and Bob agree on a finite cyclic group $G$ and a generating element $g$ in $G$.

\item Alice (the receiver) picks a random natural number $a$ and publishes $c=g^a$.

\item Bob (the sender), who wants to send a message $m \in G$ (called a
``plaintext" in cryptographic lingo) to Alice, picks a random
natural number $b$ and sends two elements, $m \cdot c^b$ and $g^b$,
to Alice. Note that $c^b=g^{ab}$.

\item   Alice recovers $m= (m \cdot c^b) \cdot ((g^b)^a)^{-1}$.

\end{enumerate}

A notable feature of the ElGamal encryption is that it is
 {\it probabilistic}, meaning
that a single plaintext can be encrypted to many possible
ciphertexts.

We also point out that the ElGamal encryption has an average
{\it expansion factor} of 2, meaning that the ciphertext is about twice as large as the corresponding plaintext.

\section{The conjugacy problem}
\label{conjugacy}

Let $G$ be a group with  solvable word problem. For $w, a \in G$,
the notation $w^a$ stands for $a^{-1}wa$. Recall that the {\it
conjugacy problem} (or {\it conjugacy decision problem}) for $G$ is:
given two elements $u, v \in G$, find out whether there is $x \in G$
such that $u^x=v$. On the other hand, the  {\it conjugacy search
problem} (sometimes also called the {\it conjugacy witness problem})
is: given two elements $a, b \in G$   and the information that
$u^x=v$ for some $x \in G$, find at least one particular element $x$
like that.

The conjugacy decision  problem is of great interest in group
theory. In contrast, the  conjugacy search problem is of interest in
complexity theory, but of little interest in group theory. Indeed,
if you know that $u$ is conjugate to $v$, you can just go over words
of the form $u^x$ and compare them to $v$ one at a time, until you
get a match. (We implicitly use here an obvious fact that a group
with solvable conjugacy problem also has solvable word problem.)
This straightforward  algorithm is at least exponential-time in the
length of $v$, and  therefore is considered infeasible for practical
purposes.

Thus, if no other algorithm is known for the  conjugacy search
problem in a group $G$, it is not unreasonable to claim that $x \to
u^x$ is a one-way function and try to build a (public-key)
cryptographic protocol on that. In other words, the assumption here
would be that in some groups $G$, the following problem is
computationally hard: given two elements $a, b$ of  $G$ and the
information that $a^x=b$ for some $x \in G$, find at least one
particular element $x$ like that. The (alleged) computational
hardness of this  problem in some particular groups (namely, in
braid groups) has been  used in several group based cryptosystems,
most notably in \cite{AAG}  and \cite{KLCHKP}. However, after some
initial excitement (which has even resulted in naming a new area of
``braid group cryptography", see e.g. \cite{Dehornoy}), it seems now
that the conjugacy search problem in a braid group may not provide
sufficient level of security; see e.g. \cite{HS, practical, MSUPKC}
for various attacks.

We start with a simple key exchange protocol, due to Ko, Lee et al.
\cite{KLCHKP}, which is modeled on the Diffie-Hellman key exchange protocol, see Section \ref{DH}.
\begin{enumerate}

 \item   An  element $w \in G$ is published.

 \item    Alice picks a private $a \in G$ and sends $w^a$
to Bob.

 \item   Bob picks a private $b \in G$ and sends $w^b$
to Alice.

 \item   Alice computes $K_A=(w^b)^a = w^{ba}$, and Bob
computes $K_B=(w^a)^b = w^{ab}$.
\end{enumerate}

If $a$ and $b$ are chosen from a pool of commuting elements of the
group $G$, then $ab=ba$, and therefore, Alice and Bob get a common
private key $K_B=w^{ab}=w^{ba}=K_A$. Typically, there are two public
subgroups $A$ and $B$ of the group $G$, given by their (finite)
generating sets, such that  $ab=ba$ for any  $a \in A$, $b \in B$.

In the paper \cite{KLCHKP}, the platform group $G$ was the braid
group $B_n$ which has some natural commuting subgroups. Selecting a
suitable platform group for the above protocol is a very nontrivial
matter; some requirements on such a group were put forward in
\cite{Shpilrain}:

\begin{enumerate}

\item[(P0)]  The conjugacy
(search)   problem in the platform group either has to be well studied
or can be reduced to a well-known problem (perhaps, in some other
area of mathematics).

\item[(P1)] The word problem in $G$ should have
a fast (at most quadratic-time) solution by a deterministic
algorithm. Better yet, there should  be an efficiently computable
``normal form" for elements of $G$.

This is required for an efficient common key extraction by
legitimate parties in a key establishment protocol, or for the
verification step in an authentication protocol, etc.

\item[(P2)]  The conjugacy search   problem should {\it not}
have an efficient solution by a deterministic algorithm.

We point out here that {\it proving} a group to have  (P2) should be
extremely difficult, if not impossible.  The  property (P2) should
therefore be considered in conjunction with (P0), i.e., the only
realistic evidence of a group $G$ having the property (P2) can be
the fact that sufficiently many people have been studying the
conjugacy (search)  problem  in $G$ over a sufficiently long time.

The next property is somewhat informal, but it is of great
importance for practical implementations:

\item[(P3)]  There should be a way to disguise elements of
$G$ so that it would be  impossible to recover  $x$ from $x^{-1}wx$
just by inspection.

One way to achieve this is to have a  {\it
normal form} for elements of $G$, which usually means that there is
an algorithm that transforms any input $u_{in}$, which is a word in
the generators of $G$, to an output $u_{out}$, which is another word
in the generators of $G$, such that $u_{in}=u_{out}$ in the group
$G$, but this is hard to detect by inspection.

In the absence of a normal form, say if $G$ is just given by means
of generators and relators without any additional information about
properties of $G$, then at least some of these relators should be
very short to be used in a disguising procedure.

\end{enumerate}

To this one can add that the platform group should not have a linear
representation of a small dimension since otherwise, a linear
algebra attack might be feasible.

\subsection{The Anshel-Anshel-Goldfeld key exchange protocol}
\label{AAG}

 In this section, we are going to describe a key establishment
 protocol from \cite{AAG} that really  stands out because, unlike other protocols
 based on the (alleged) hardness of the conjugacy search problem,
it does not employ any commuting or commutative subgroups of a given
 platform group and can, in fact, use any non-abelian
group with efficiently  solvable word problem as the platform. This
really makes a difference and gives a big advantage to the protocol
of \cite{AAG} over most protocols in this and the following section.
The choice of the platform group $G$ for this protocol is a delicate
matter though. In the original paper \cite{AAG}, a braid group was
suggested as the platform, but with this platform the protocol was
subsequently attacked in several different ways, see e.g.
\cite{Ben-Zvi}, \cite{GKTTV2}, \cite{GKTTV}, \cite{HS},
\cite{LeePark}, \cite{LeeLee1}, \cite{practical}, \cite{MSUPKC},
\cite{Tsaban}. The search for a good platform group for this
protocol still continues.

Now we give a description of the AAG protocol.  A group $G$ and
elements $a_1,...,a_k, b_1,...,b_m \in G$ are public.

\begin{enumerate}
 \item[{\bf (1)}] Alice picks a private $x \in G$ as a word in
$a_1,...,a_k$ (i.e.,  $x=x(a_1,...,a_k)$)  and sends
$b_1^x,...,b_m^x$ to Bob.

 \item[{\bf (2)}] Bob picks a private $y \in G$
as a word in $b_1,...,b_m$ and sends $a_1^y,...,a_k^y$ to Alice.

 \item[{\bf (3)}] Alice computes $x(a_1^y,...,a_k^y) = x^y = y^{-1}xy$,
and Bob computes \\
$y(b_1^x,...,b_m^x) = y^x = x^{-1}yx$. Alice and Bob then come up
with a common private key $K=x^{-1}y^{-1}xy$ (called the
\index{ind}{Commutator} \emph{commutator} of $x$ and $y$) as
follows: Alice multiplies $y^{-1}xy$ by $x^{-1}$ on the left,  while
Bob multiplies $x^{-1}yx$ by $y^{-1}$ on the left, and then takes
the inverse of the whole thing: $(y^{-1}x^{-1}yx)^{-1} =
x^{-1}y^{-1}xy$.

\end{enumerate}

It may seem that solving the (simultaneous) conjugacy search problem
for $b_1^x,...,b_m^x; a_1^y,...,a_k^y$ in the group $G$ would allow
  an adversary
to get the secret key $K$. However,  if we look at Step (3) of the
protocol, we see that the adversary would have to know either $x$ or
$y$ not simply as  a word in the generators of the group $G$, but as
a word in $a_1,...,a_k$ (respectively, as  a word in $b_1,...,b_m$);
otherwise, he  would not be able to compose, say, $x^y$ out of
$a_1^y,...,a_k^y$. That means the adversary would also have to solve
the  {\it membership search
problem}:

\begin{quote}
 Given   elements $x, a_1,...,a_k$ of a group $G$,
find an expression (if it exists) of  $x$ as a word in
$a_1,...,a_k$.
\end{quote}

 We note  that the   membership {\it  decision} problem is to
determine whether or not a given $x \in G$ belongs to the subgroup
of $G$ generated by given $a_1,...,a_k$. This problem turns out to
be quite hard in many groups. For instance, the membership decision
problem in a braid group $B_n$ is algorithmically unsolvable if $n
\ge 6$ because such a braid group contains
 subgroups isomorphic to  $F_2 \times F_2$ (that would
be, for example, the subgroup generated by $\sigma_1^2, \sigma_2^2,
\sigma_4^2$, and  $\sigma_5^2$, see \cite{Collins}), where $F_2$ is
the free group of rank 2. In the group $F_2 \times F_2$, the
membership decision problem is algorithmically unsolvable by an old
result of Mihailova \cite{Mihailova}.

 We also note that if the adversary finds, say, some $x' \in G$
such that ~$b_1^x=b_1^{x'}, ... ,  b_m^x=b_m^{x'}$, there is no
guarantee that $x'=x$ in $G$.  Indeed, if $x'=c_b x$, where $c_b
b_i= b_i c_b$ for all $i$ (in which case we say that $c_b$ {\it
centralizes} $b_i$), then $b_i^x=b_i^{x'}$ for all $i$, and
therefore  $b^x=b^{x'}$ for any element $b$  from the subgroup
generated by $b_1,...,b_m$; in particular, $y^x=y^{x'}$. Now the
problem is that if $x'$ (and, similarly, $y'$)
 does not belong to the
subgroup $A$ generated by $a_1,...,a_k$ (respectively, to the
subgroup $B$ generated by $b_1,...,b_m$),
 then the adversary   may not obtain the correct common secret key $K$.
On the other hand, if  $x'$ (and, similarly, $y'$)  does   belong to
the subgroup $A$ (respectively, to the subgroup $B$), then the
adversary   will  be able to get the correct $K$ even though his
$x'$ and $y'$ may be different from $x$ and $y$, respectively.
Indeed, if $x'=c_b x$, $y'=c_a y$, where $c_b$ centralizes $B$  and
 $c_a$ centralizes $A$ (elementwise), then
$$(x')^{-1}(y')^{-1}x'y'=(c_b x)^{-1}(c_a y)^{-1}c_b x c_a y=x^{-1}c_b^{-1}y^{-1}c_a^{-1}c_b x c_a y=
x^{-1}y^{-1}xy=K$$
 because $c_b$ commutes with $y$ and with  $c_a$ (note that $c_a$ belongs
to the subgroup $B$, which follows from the assumption $y'=c_a y \in
B$, and, similarly, $c_b$ belongs to $A$), and $c_a$ commutes with
$x$.

We emphasize that the adversary ends up with the corrrect key $K$
(i.e., $K=(x')^{-1}(y')^{-1}x'y'=x^{-1}y^{-1}xy$) {\it if and only
if} $c_b$ commutes with $c_a$. The only visible way to ensure this
is to have $x' \in A$ and  $y' \in B$. Without verifying at least
one of these inclusions, there seems to be no way for the adversary
to make sure that he got the correct key.

Therefore, it appears that if the adversary chooses to solve the
conjugacy search problem in the group $G$ to recover $x$ and $y$, he
will then have to face either the  membership search problem or
  the  membership decision problem; the  latter  may very well be
algorithmically unsolvable in a given group. The bottom line is that
the adversary should actually be solving a (probably) more difficult
(``subgroup-restricted") version of the conjugacy search problem:

\begin{quote}
 Given a group $G$, a subgroup  $A \le G$, and two
 elements $g, h \in G$,  find $x \in A$
such that $h=x^{-1}gx$, given that at least one such $x$ exists.
\end{quote}

\subsection{The  twisted conjugacy problem}
\label{twist}

Let $\phi, \psi$ be two fixed automorphisms (more generally,
endomorphisms) of a group $G$. Two elements $u, v \in G$  are called
{\em $(\phi,\psi)$-double-twisted conjugate} if there is an element
$w \in G$ such that  $uw^\phi = w^\psi v$.  When $\psi = id$, then
$u$ and $v$ are called {\em $\phi$-twisted conjugate}, while in the
case $\phi = \psi = id$,  $u$ and $v$ are just usual conjugates of
each other.

The twisted (or double twisted) conjugacy problem in $G$ is:

\begin{quote}  decide whether or not two  given elements $u, v \in G$ are
twisted (double twisted) conjugate  in $G$ for a fixed pair of
endomorphisms $\phi, \psi$ of the group $G$.
\end{quote}

Note that if $\psi$ is an automorphism, then
$(\phi,\psi)$-double-twisted conjugacy problem reduces to
$\phi\psi^{-1}$-twisted conjugacy problem, so in this case  it is
sufficient to consider just the twisted conjugacy problem. This
problem was studied from group-theoretic perspective, see e.g.
\cite{Ventura-Romankov:2009,  Romankov:2011,
Fel'shtyn-Leonov-Troitsky}, and in  \cite{twisted conjugacy} it was
used in an authentication protocol. It is interesting that the
research in \cite{Ventura-Romankov:2009, Romankov:2011} was probably
motivated by cryptographic applications, while the authors of
\cite{Fel'shtyn-Leonov-Troitsky} arrived at the twisted conjugacy
problem motivated by problems in topology.

\section{The decomposition problem}
\label{decomposition}

Another ramification of the conjugacy search  problem is
the following  {\it decomposition search problem}:

\begin{quote}
Given  two elements $w$ and $w'$ of a group $G$, find two  elements
$x \in A$ and $y \in B$ that would belong to  given subsets (usually
subgroups) $A, B \subseteq G$ and satisfy $x\cdot w\cdot y = w'$,
provided at least one such pair of elements exists.
\end{quote}

We  note that if in the above problem $A=B$ is a subgroup, then this
problem is also known as the {\it double coset problem}.

We also  note that {\it some} $x$ and $y$ satisfying the equality
$x\cdot w\cdot y =   w'$ always exist (e.g.  $x=1, ~y=w^{-1}w'$), so
the point is to have them satisfy the conditions $x \in A$ and $y
\in B$. We therefore will not usually refer to this problem as a
{\it subgroup-restricted} decomposition search problem because it is
always going to be subgroup-restricted; otherwise it does not make
much sense. We also note that the most commonly considered special
case of the decomposition search problem so far is where $A=B$.

We are going to show in Section \ref{relations} that solving the
conjugacy search problem is unnecessary for an adversary to get the
common secret key in the Ko-Lee (or any similar) protocol (see our
Section \ref{conjugacy}); it is sufficient to solve a seemingly
easier decomposition search problem. This was mentioned, in passing,
in the paper \cite{KLCHKP}, but the significance of this observation
was downplayed there.

We note that the membership condition $x, y \in A$ may not be easy to verify
for some subsets $A$. The authors of \cite{KLCHKP} do not
address this problem; instead they mention, in justice, that if one
uses a ``brute force" attack by simply going over elements of $A$
one at a time, the above condition  will be satisfied automatically.
This however may not be the case with other, more practical,
attacks.

We also note that the conjugacy search problem is a special case of
the decomposition problem where $w'$ is conjugate to $w$ and
$x=y^{-1}$. The claim that the decomposition problem should be
easier than the conjugacy search  problem is intuitively clear since
it is generally  easier to solve an equation with two unknowns than
a special case of the same equation with just one unknown. We admit
however that there might be exceptions to this general rule.

Now we give  a formal description of  a typical  protocol based on
the decomposition problem. There is a public group $G$, a public
element $w \in G$, and two public subgroups $A, B \subseteq G$
commuting elementwise, i.e., $ab=ba$ for any $a \in A, b \in B$.

\begin{enumerate}

 \item Alice randomly selects private elements $a_1, a_2
\in A$. Then she sends  the  element $a_1 w a_2$  to Bob.

 \item  Bob randomly selects private elements  $b_1, b_2 \in
B$. Then he sends the element $b_1 w b_2$  to Alice.

 \item   Alice computes $K_A= a_1 b_1 w b_2 a_2$, and Bob computes $K_B= b_2 a_1 w b_1 a_2$. Since
$a_ib_i=b_ia_i$ in $G$, one has $K_A=K_B=K$ (as an element of $G$),
which is now Alice's and Bob's common secret key.

\end{enumerate}

We now discuss several modifications of the above protocol.

\subsection{``Twisted"  protocol}
\label{twist_protocol}

This idea is due to Shpilrain and  Ushakov \cite{decomposition}; the
following modification of the above protocol appears to be more
secure (at least for some choices of the platform group) against
so-called ``length based" attacks (see e.g. \cite{GKTTV2},
\cite{GKTTV}, \cite{HS}), according to computer experiments. Again,
there is a public group $G$ and two public subgroups $A, B \leq G$
commuting elementwise.

\begin{enumerate}

\item Alice randomly selects private elements $a_1 \in A$ and $b_1 \in
B$. Then she sends  the  element $a_1 w b_1$ to  Bob.

\item Bob randomly selects private elements  $b_2 \in B$ and $a_2 \in
A$. Then he sends  the  element $b_2 w a_2$  to  Alice.

\item Alice computes $K_A= a_1 b_2 w a_2 b_1 = b_2 a_1 w b_1 a_2$, and Bob
computes $K_B= b_2 a_1 w b_1 a_2$. Since $a_ib_i=b_ia_i$ in $G$, one
has $K_A=K_B=K$ (as an element of $G$), which is now Alice's and
Bob's common secret key.

\end{enumerate}

\subsection{Finding intersection of given subgroups}
\label{central}

Another modification of the protocol in Section \ref{decomposition}
is also due to Shpilrain and  Ushakov \cite{decomposition}. First we
give a sketch of the idea.

Let $G$ be a group and $g \in G$. Denote by $C_G(g)$  the {\em
centralizer}  of $g$ in $G$, i.e., the set of elements $h \in G$
such that $hg = gh$. For $S = \{g_1,\ldots,g_k \} \subseteq G$,
$C_G(g_1,\ldots,g_k)$ denotes the centralizer of $S$ in $G$, which
is the intersection of  the centralizers $C_G(g_i), i=1,...,k$.

Now, given a public $w \in G$, Alice
privately selects $a_1 \in G$ and publishes a subgroup   $B
\subseteq C_G(a_1)$ (we tacitly assume here that $B$ can be computed
efficiently). Similarly, Bob privately selects $b_2 \in G$ and
publishes a subgroup $A \subseteq C_G(b_2)$. Alice then selects $a_2
\in A$ and sends $w_1=a_1wa_2$ to Bob, while Bob selects $b_1 \in B$
and sends $w_2=b_1wb_2$ to Alice.

Thus,  in the first  transmission, say, the adversary faces the
problem of finding $a_1, a_2$ such that $w_1 = a_1 w a_2$, where
$a_2 \in A$, but there is no explicit indication of where to choose
$a_1$ from. Therefore, before arranging something like a length
based attack in this case, the adversary would have to compute
generators of the centralizer $C_G(B)$ first (because $a_1 \in
C_G(B)$), which is usually a hard problem by itself since it
basically amounts to finding the intersection of  the centralizers
of individual elements, and finding (the generators of) the
intersection of subgroups is a  notoriously difficult problem for
most groups considered in combinatorial group theory.

Now we give a formal description of the protocol from
\cite{decomposition}. As usual, there is a public group $G$, and let
$w \in G$ be public, too.

\begin{enumerate}

\item Alice chooses an element $a_1 \in G$, chooses a
subgroup of $C_G(a_1)$, and publishes its generators $A = \{
\alpha_1,\ldots,\alpha_k\}$.

\item Bob chooses an element $b_2 \in G$, chooses a subgroup
of $C_G(b_2)$, and publishes  its generators $B =
\{\beta_1,\ldots,\beta_m \}$.

\item Alice chooses a random element $a_2$ from
$gp<{\beta_1,\ldots,\beta_m}>$ and sends   $P_A =a_1 w a_2$ to Bob.

\item Bob chooses a random element $b_1$ from
$gp<{\alpha_1,\ldots,\alpha_k}>$ and sends  $P_B = b_1 w b_2$ to
Alice.

\item   Alice computes $K_A = a_1 P_B a_2$.

\item  Bob computes $K_B = b_1 P_A b_2$.

\end{enumerate}

Since $a_1b_1 = b_1a_1$ and $a_2b_2 = b_2a_2$, we have $K = K_A =
K_B$, the shared secret  key.

We note that in \cite{Tsaban}, an attack on this protocol was
offered (in the case where a braid group is used as the platform),
using what the author calls the {\it linear centralizer method}.
Then, in \cite{Ben-Zvi}, another method of cryptanalysis (called the
{\it algebraic span cryptanalysis}) was offered, applicable to
platform groups that admit an efficient linear representation. This
method yields attacks on various protocols, including the one in
this section, if a braid group is used as the platform.

\subsection{Commutative subgroups}
\label{Commuative_subgroups}

Instead of using {\it commuting} subgroups $A, B \le G$, one can use
{\it commutative} subgroups. Thus, suppose $A, B \le G$ are two
public commutative subgroups (or subsemigroups) of a group $G$, and
let $w \in G$ be a public element.

\begin{enumerate}

\item Alice randomly selects private elements $a_1 \in A$, $b_1 \in B$.
 Then she sends  the  element $a_1 w b_1$  to Bob.

\item  Bob randomly selects private elements $a_2 \in A$, $b_2 \in B$.
Then he sends the element $a_2 w b_2$  to Alice.

 \item   Alice computes $K_A= a_1 a_2 w b_2 b_1$, and Bob computes $K_B= a_2 a_1 w b_1 b_2$. Since
$a_1 a_2=a_2 a_1$ and $b_1 b_2=b_2 b_1$  in $G$, one has $K_A=K_B=K$
(as an element of $G$), which is now Alice's and Bob's common secret
key.

\end{enumerate}

\subsection{The factorization problem}
\label{factorization}

The  {\it factorization search problem} is a special case of the decomposition search  problem:

\begin{quote}
Given  an element $w$ of a group $G$ and two subgroups
$A, B \leq G$, find any two elements $a \in A$ and $b \in B$ that
would   satisfy $a\cdot  b = w$, provided at least one such pair of elements exists.
\end{quote}

The following protocol relies in its security on the computational
hardness of the factorization search problem.  As before, there is a
public group $G$, and two public subgroups $A, B \leq G$ commuting
elementwise, i.e., $ab=ba$ for any $a \in A, b \in B$.

\begin{enumerate}

 \item Alice randomly selects private elements $a_1 \in A, b_1
\in B$. Then she sends  the  element $a_1 b_1$  to Bob.

 \item  Bob randomly selects private elements  $a_2 \in A, b_2 \in
B$. Then he sends the element $a_2 b_2$  to Alice.

 \item   Alice computes
$$K_A=  b_1 (a_2 b_2) a_1 = a_2 b_1 a_1 b_2 = a_2 a_1 b_1 b_2,$$
and Bob computes
$$K_B=  a_2 (a_1 b_1) b_2 = a_2 a_1 b_1 b_2.$$
Thus, $K_A=K_B=K$  is now Alice's and Bob's common secret key.
\end{enumerate}

We note that the adversary, Eve, who knows the elements  $a_1 b_1$
and $a_2 b_2$, can compute  $(a_1 b_1)(a_2 b_2)=a_1 b_1 a_2 b_2= a_1
a_2 b_1 b_2$ and $(a_2 b_2)(a_1 b_1) = a_2 a_1 b_2 b_1$, but neither
of these products is equal to $K$ if $a_1 a_2 \ne a_2 a_1$ and $b_1
b_2 \ne b_2 b_1$.

Finally, we point out a {\it decision} factorization problem:

\begin{quote}
 Given  an element $w$ of a group $G$ and two subgroups
$A, B \leq G$, find out whether or not there are  two elements $a
\in A$ and $b \in B$  such that $w = a\cdot  b$.
\end{quote}

This seems to be a new and  non-trivial algorithmic problem in
group theory, motivated by cryptography.

\section{The word problem}
\label{WP}

The word problem ``needs no introduction", but it probably makes
sense to spell out the {\it word search problem:}

\begin{quote} Suppose $H$ is a group given by a finite presentation $<X; R>$ and
let $F(X)$ be the free group with the set $X$ of free generators.
Given a group word $w$ in the alphabet $X$, find a sequence of
conjugates of elements from $R$ whose product is equal to $w$ in the
free group $F(X)$.
\end{quote}

Long time ago, there was an attempt to use the undecidability of the
{\it decision} word problem (in some groups) in public key
cryptography \cite{MW}. This was, in fact, historically the first
attempt to employ a hard algorithmic  problem from combinatorial
group theory in public key cryptography. However, as was pointed out
in \cite{BMS}, the problem that is actually used in \cite{MW} is not
the word problem, but the {\it  word  choice problem}: given  $g,
w_1, w_2 \in G$, find out whether $g =w_1$  or   $g =w_2$ in $G$,
provided one of the two equalities holds. In this problem, both
parts are recursively solvable for any recursively presented
platform group $G$ because they both are the ``yes" parts of the
word problem. Therefore, undecidability of the actual word problem
in the platform group has no bearing on the security of the
encryption scheme in \cite{MW}.

On the other hand, employing decision problems (as opposed to search
 problems) in public-key cryptography would allow one to depart
from the canonical paradigm and construct cryptographic protocols
with new properties, impossible in the canonical model. In
particular, such protocols can be secure against some ``brute force"
attacks by a computationally unbounded adversary. There is a price
to pay for that, but the price is reasonable: a legitimate receiver
decrypts correctly with probability that can be made very close to
1, but not equal to 1. This idea was implemented in \cite{Osin}, so
the exposition below follows that paper.

We assume that the sender (Bob) is given a presentation $\Gamma$
(published by the receiver Alice) of a group $G$ by generators and
defining relators:
$$\Gamma = \langle x_1, x_2,\ldots, x_n  \mid r_1, r_2,\dots \rangle.$$

No further information about the group $G$ is available to Bob.

Bob is instructed to transmit his private bit to Alice by
transmitting a word $u=u(x_1, \dots, x_n)$ equal to 1 in $G$  in
place of ``1" and a word  $v=v(x_1, \dots, x_n)$  not equal to 1 in
$G$ in place of ``0".

Now we have to specify the  algorithms that Bob should use to select
his words.
\medskip

\noindent  {\bf Algorithm ``0"} (for selecting a word  $v=v(x_1,
\dots, x_n)$  not equal to 1 in $G$)  is quite simple: Bob just
selects a random word by building it letter-by-letter, selecting
each letter uniformly from the set $X = \{x_1, \dots, x_n, x_1^{-1},
\dots, x_n^{-1}\}$. The length of such a word  should be a random
integer from an interval that Bob selects up front, based on his
computational abilities.
\medskip

\noindent  {\bf Algorithm ``1"} (for selecting a word  $u=u(x_1,
\dots, x_n)$  equal to 1 in $G$) is slightly more complex. It
amounts to applying a random sequence of operations of the following
two kinds, starting with the empty word:

\begin{enumerate}

    \item Inserting into a random place in the current word  a pair
$hh^{-1}$ for a random word  $h$.

\item    Inserting into a random place in the current word  a random conjugate
    $g^{-1}r_ig$ of a random defining relator  $r_i$.

\end{enumerate}

The length of the resulting word  should be in the same range as the
length of the output of Algorithm ``0", for indistinguishability.

\subsection{Encryption emulation attack} Now let us see what happens
if a  computationally unbounded adversary uses what is called {\it
encryption emulation attack} on Bob's encryption. This kind of
attack always succeeds against ``traditional" encryption protocols
where the receiver decrypts correctly with probability exactly 1.
The encryption emulation attack is:

\begin{quote}
For either bit,  generate its encryption over and over again, each
time with fresh randomness, until the ciphertext to be attacked is
obtained. Then the corresponding plaintext is the bit that was
encrypted.
\end{quote}

Thus, the (computationally unbounded) adversary is building up two
lists, corresponding to two algorithms above. Our first observation
is that the list that corresponds to the Algorithm ``0" is useless
to the adversary because it is eventually going to contain {\it all}
words in the alphabet $X = \{x_1, \dots, x_n, x_1^{-1}, \dots,
x_n^{-1}\}$. Therefore, the adversary may just as well forget about
this list and focus on the other one, that corresponds to the
Algorithm ``1".

Now the situation boils down to the following: if a  word  $w$
 transmitted by Bob appears on the list, then it is equal to 1 in
 $G$. If not, then not. The only problem is: how can one
 conclude that $w$ does {\it not} appear on the list if the list is
 infinite?  Of course, there is no
 infinity in real life, so the list is actually finite because of
 Bob's computational limitations. Still, at least in theory,
the adversary does not know a {\it bound} on the size of the list if
she does not know Bob's computational limits.

Then, perhaps the adversary can stop at some point and conclude
that $w \ne 1$ with overwhelming probability, just like Alice does?
The point however is  that this probability may not at all be as
``overwhelming" as the probability of the  correct decryption by
Alice. Compare:

\begin{enumerate}

    \item For Alice to decrypt correctly ``with overwhelming
probability", the probability $P_1(N)$
    for a random word  $w$ of length  $N$ not to be equal to 1
    should converge to 1 (reasonably fast) as $N$ goes to infinity.

\item For the adversary  to decrypt correctly ``with overwhelming
probability", the probability $P_2(N, f(N))$ for a random word  $w$
of length  $N$ produced by the Algorithm ``1" to have a {\it  proof}
of length $\le f(N)$ verifying that $w=1$, should converge to 1
 as $N$ goes to infinity. Here $f(N)$ represents
the adversary's computational capabilities; this function can be
arbitrary, but fixed.

\end{enumerate}

We  see that the functions $P_1(N)$ and  $P_2(N)$ are of very
 different nature, and any correlation between them is unlikely.
We note that the function $P_1(N)$ is generally well understood, and
in particular, it is known that in any infinite group $G$, $P_1(N)$
indeed converges to 1 as $N$ goes to infinity.

On the other hand, functions $P_2(N, f(N))$ are more complex; note
 also that they may depend on a particular algorithm used by Bob to
 produce words equal to 1. The Algorithm ``1" described in this
 section is very straightforward; there are more delicate algorithms
 discussed in \cite{MSUbook2}.

Functions $P_2(N, f(N))$  are currently subject of active research,
and  in particular, it appears likely that there are groups in which
$P_2(N, f(N))$ does not converge to 1 at all, if an algorithm used
to produce words equal to 1 is chosen intelligently.

We also note in passing that if in a group $G$ the word problem is
recursively unsolvable, then the length of a proof verifying that
$w=1$ in $G$ is  not bounded by any recursive function of the length
of $w$.

 Of course, in real life, the adversary may know a bound on the size of the
 list based on a general idea of what kind of hardware
 may be available to Bob; but then again, in real life the adversary would be
 computationally bounded, too. Here we note (again, in passing)
that there are groups $G$ with efficiently solvable word problem and
words $w$ of length  $n$  equal to 1 in $G$, such that the length of
a proof verifying that $w=1$ in $G$ is  not bounded by any tower of
exponents in $n$, see \cite{P}.

Thus, the bottom line is: in theory, the adversary cannot positively
identify the bit that Bob has encrypted by a word  $w$ if she just
uses the ``encryption emulation" attack.
 In fact, such an identification would be equivalent to solving the {\it word
 problem} in $G$, which would contradict the well-known fact that
 there are (finitely presented) groups with recursively unsolvable word
 problem.

It would be nice, of course,  if the adversary was unable to
positively decrypt using ``encryption  emulation" attacks even if
she {\it did} know Bob's computational limitations. This, too, can
be arranged, see the following subsection.

\subsection{Encryption: trick and treat} \label{encryption}

Building on the ideas from the previous subsection and combining
them with a simple yet subtle trick, we describe here an encryption
protocol from \cite{Osin} that has the following features:

 \begin{itemize}
\item[(F1)] Bob encrypts his private bit sequence by a word in a public
alphabet $X$.

\item[(F2)] Alice (the receiver)   decrypts Bob's transmission  correctly with probability
that can be made arbitrarily close to 1, but not equal to 1.

\item[(F3)] The adversary, Eve, is assumed to have no bound on the speed of
computation or on the storage space.

\item[(F4)] Eve is assumed to have complete information on the algorithm(s) and
hardware that Bob uses for encryption. However, Eve cannot predict
outputs of Bob's random numbers generator.

\item[(F5)] Eve cannot decrypt Bob's bit  correctly
with probability $>\frac{3}{4}$  by emulating Bob's encryption
algorithm.

\end{itemize}

This leaves Eve with the only possibility:  to attack Alice's
decryption algorithm or her algorithm for obtaining public keys, but
this is a different story. Here we only discuss the encryption
emulation attack, to make a point that this attack can be
unsuccessful if the probability of the legitimate decryption is
close to 1, but not exactly 1.

Here is the relevant protocol (for encrypting a single bit).

 \begin{itemize}
\item[(P0)]  Alice publishes two presentations:

$$\Gamma_1 = \langle x_1, x_2,\ldots, x_n  \mid r_1, r_2,\dots \rangle$$
\vskip -.2cm

$$\Gamma_2 = \langle x_1, x_2,\ldots, x_n  \mid s_1, s_2,\dots \rangle.$$

\noindent  One of them defines the trivial group, whereas the other
one defines an infinite group, but only Alice knows which one is
which. Bob is instructed to transmit his private bit to Alice as
follows:

\item[(P1)] In place of ``1", Bob transmits a pair of words $(w_1, w_2)$ in the
alphabet\\
 $X=\{x_1, x_2,\ldots, x_n, x_1^{-1}, \ldots, x_n^{-1}\}$, where
$w_1$ is selected randomly, while $w_2$ is selected to be equal to 1
in the group $G_2$ defined by $\Gamma_2$ (see e.g. Algorithm ``1" in the
previous section).

\item[(P2)] In place of ``0", Bob transmits a pair of words $(w_1,
w_2)$, where $w_2$ is selected randomly, while $w_1$ is selected to
be  equal to 1 in the group $G_1$ defined by $\Gamma_1$.

\end{itemize}

Under our assumptions (F3), (F4) Eve can identify the word(s) in the
transmitted pair which is/are equal to 1 in the corresponding
presentation(s), as well as the word, if any, which is not equal to
1. There are the following possibilities:

 \begin{enumerate}

  \item $w_1=1$ in $G_1$, $w_2=1$ in $G_2$;

\item $w_1=1$ in $G_1$, $w_2 \ne 1$ in $G_2$;

\item $w_1\ne 1$ in $G_1$, $w_2 = 1$ in $G_2$.

\end{enumerate}

 It is easy to see that the possibility (1) occurs with probability
 $\frac{1}{2}$ (when Bob wants to transmit ``1" and $G_1$ is
 trivial, or when Bob wants to transmit ``0" and $G_2$ is
 trivial). If this possibility occurs, Eve cannot decrypt Bob's bit correctly
with probability $>\frac{1}{2}$. Indeed, the only way for Eve  to
decrypt in this case would be to find out which presentation
$\Gamma_i$ defines the trivial group, i.e., she would have to attack
Alice's algorithm for obtaining a public key, which would not be
part of the encryption emulation attack anymore. Here we just  note,
in passing, that there are many different ways to construct
presentations of the trivial group, some of them involving a lot of
random choices. See e.g. \cite{MMS} for a survey on the subject.

In any case, our claim (F5) was that Eve cannot decrypt Bob's bit
correctly with probability $>\frac{3}{4}$  by emulating Bob's
encryption algorithm, which is obviously true in this scheme since
the probability for Eve to decrypt correctly is, in fact, precisely
$\frac{1}{2} \cdot \frac{1}{2} + \frac{1}{2}  \cdot 1 =
\frac{3}{4}$. (Note that Eve  decrypts correctly with probability 1
if either of the possibilities (2) or (3) above occurs.)

\section{The isomorphism inversion problem}
\label{isomorphism}

The isomorphism (decision) problem for groups is very well known: suppose two
groups are given by their finite presentations in terms of
generators and defining relators, then find out whether  the groups are
isomorphic. The search version of this problem is well known, too:
given two finite presentations defining isomorphic groups, find a
particular isomorphism between the groups.

Now the following problem, of interest in cryptography, is not what
was previously considered in combinatorial group theory:

\begin{quote} Given two finite presentations defining isomorphic
groups, $G$  and $H$,  and an isomorphism $\varphi: G \to H$, find
$\varphi^{-1}$.
\end{quote}


Now we describe an encryption scheme whose security is based on the
alleged computational hardness of the isomorphism inversion problem.
Our idea itself is quite simple: encrypt with a public isomorphism
$\varphi$ that is computationally infeasible for the adversary to
invert. A legitimate receiver, on the other hand, can efficiently
compute $\varphi^{-1}$ because she knows a factorization of
$\varphi$ in a product of ``elementary'', easily invertible,
isomorphisms.

What is interesting to note is that this encryption is {\it
homomorphic} because $\varphi(g_1g_2)=\varphi(g_1)\varphi(g_2)$ for
any $g_1, g_2 \in G$. The significance of this observation is due to
a result of \cite{OS}: if the group $G$ is a non-abelian finite simple
group, then any homomorphic encryption on $G$ can be converted to a
{\it fully homomorphic encryption} (FHE) scheme, i.e., encryption
that respects not just one but two operations: either boolean AND
and OR or arithmetic addition and multiplication.

In summary, a relevant scheme can be built as follows. Given a
public presentation of a group $G$ by generators and defining
relations, the receiver (Alice) uses  a chain of private
``elementary'' isomorphisms $G \to H_1 \to ... \to H_k \to H$, each
of which is easily invertible, but the (public) composite
isomorphism $\varphi : G \to H$ is hard to invert without the
knowledge of a factorization in a product of ``elementary'' ones.
(Note that $\varphi$ is published as a map taking the generators of
$G$ to words in the generators of $H$.) Having obtained this way a
(private) presentation $H$, Alice discards some of the defining
relations to obtain a public presentation $\hat H$. Thus, the group
$H$, as well as the group $G$ (which is isomorphic to $H$), is a
homomorphic image of the group $\hat H$. (Note that $\hat H$ has the
same set of generators as $H$ does but has fewer defining
relations.) Now the sender (Bob), who wants to encrypt his plaintext
$g \in G$, selects an arbitrary word $w_g$ (in the generators of
$G$) representing the element $g$ and applies the public isomorphism
$\varphi$ to $w_g$ to get $\varphi(w_g)$, which is a word in the
generators of $H$ (or $\hat H$, since $H$ and $\hat H$ have the same
set of generators). He then selects an arbitrary word $h_g$ in the
generators of $\hat H$ representing the same element of $\hat H$ as
$\varphi(w_g)$ does, and this is now his ciphertext:  $h_g = E(g)$.
To decrypt, Alice applies her private map $\varphi^{-1}$ (which is a
map taking the generators of $\hat H$ to words in the generators of
$G$) to $h_g$ to get a word $w'_g = \varphi^{-1}(h_g)$. This word
$w'_g$ represents the same element of $G$ as $w_g$ does because
$\varphi^{-1}(h_g) = \varphi^{-1}(\varphi(w_g)) = w_g$ in the group
$G$ since both $\varphi$  and $\varphi^{-1}$ are homomorphisms, and
the composition of $\varphi$  and $\varphi^{-1}$ is the identity map
on the group $G$, i.e., it takes every word in the generators of $G$
to a  word representing the same element of $G$. Thus, Alice
decrypts correctly.

We emphasize here that a plaintext is a {\it group element} $g \in
G$, not a word in the generators of $G$. This implies, in
particular, that there should be some kind of canonical way (a ``normal form") of
representing elements of $G$. For example, for elements of an
alternating group $A_m$ (these groups are finite  non-abelian simple
groups if $m \ge 5$), such a canonical representation can be the
standard representation by a permutation of the set $\{1, \ldots,
m\}$.

Now we are going to give more details on how one can construct a
sequence of ``elementary'' isomorphisms starting with a given
presentation of a group $G = \langle x_1, x_2,\ldots \mid r_1,
r_2,\dots\rangle$. (Here $x_1, x_2,\ldots$ are generators and $r_1,
r_2,\dots$ are defining relators). These ``elementary'' isomorphisms
are called Tietze transformations. They are universal in the sense
that they can be applied to any (semi)group presentation. Tietze
transformations are of the following  types:

\begin{description}
\item[(T1)]
\emph{Introducing a~new generator}: Replace $\langle x_1,x_2,\ldots
\mid r_1, r_2,\dots\rangle$ by
 $\langle y, x_1,x_2,\ldots \mid
ys^{-1}, r_1, r_2,\dots\rangle$, where $s=s(x_1,x_2,\dots )$ is an
arbitrary element in the generators $x_1,x_2,\dots$.

\item[(T2)]
\emph{Canceling a~generator} (this is the converse of (T1)): If we
have a~presentation of the form $\langle y, x_1,x_2,\ldots \mid q,
r_1, r_2,\dots\rangle$, where $q$ is of the form $ys^{-1}$, and $s,
r_1, r_2,\dots $ are in the group generated by $x_1,x_2,\dots$,
replace this presentation by $\langle x_1,x_2,\ldots \mid r_1,
r_2,\dots\rangle$.

\item[(T3)]
\sloppy \emph{Applying an    automorphism}: Apply an automorphism of
the free group generated by $x_1,x_2,\dots $ to all the relators
$r_1, r_2,\dots$.

\item[(T4)]
\emph{Changing defining relators}: Replace the set $r_1, r_2,\dots$
of defining relators by another set $r_1', r_2',\dots$ with the same
normal closure. That means, each of  $r_1', r_2',\dots$ should
belong to the normal subgroup generated by $r_1, r_2,\dots$, and
vice versa.
\end{description}

Tietze proved (see e.g. \cite{LS}) that two groups given by
presentations $\langle x_1,x_2,\ldots \mid r_1, r_2,\dots\rangle$
and $\langle y_1,y_2,\ldots \mid s_1, s_2,\dots\rangle$ are
isomorphic if and only if one can get from one of the presentations
to the other by a sequence of transformations \textup{(T1)--(T4)}.

For each Tietze transformation of the types \textup{(T1)--(T3)}, it
is easy to obtain an explicit isomorphism (as a map on generators)
and its inverse. For a Tietze transformation of the type
\textup{(T4)}, the isomorphism is just the identity map. We would
like here to make Tietze transformations of the type \textup{(T4)}
recursive, because {\it a priori} it is not clear how Alice can
actually implement these transformations. Thus, Alice can use the
following recursive version of \textup{(T4)}:
\medskip

\noindent {\bf (T4$'$)} In the set $r_1, r_2,\dots$, replace some
$r_i$ by one of the:  $r_i^{-1}$,  $r_i r_j$, $r_i r_j^{-1}$, $r_j
r_i$, $r_j r_i^{-1}$, $x_k^{-1} r_i x_k$, $x_k r_i x_k^{-1}$, where
$j \ne i$, and $k$ is arbitrary.
\medskip

One particularly useful feature of Tietze transformations is that
they can break long defining relators into  short pieces (of length
3 or 4, say) at the expense of introducing more generators, as
illustrated by the following simple  example. In this example, we
start with a presentation having two relators of length 5 in 3
generators, and end up with a presentation having 4 relators of
length 3 and one relator of length 4, in 6 generators. The $\cong$
symbol below means ``is isomorphic to.''

\begin{example} \label{Tietze_example}  $G = \langle x_1, x_2, x_3 ~\mid ~x_1^2x_2^3,
~x_1x_2^2x_1^{-1}x_3 \rangle ~\cong ~\langle x_1, x_2, x_3, x_4
~\mid ~x_4=x_1^2, ~x_4x_2^3, ~x_1x_2^2x_1^{-1}x_3 \rangle\\
\cong \langle x_1, x_2, x_3, x_4, x_5 ~\mid  ~x_5=x_1x_2^2,
~x_4=x_1^2, ~x_4x_2^3, ~x_5x_1^{-1}x_3 \rangle \cong $ (now switching $x_1$ and  $x_5$ -- this is (T3))\\
$\cong \langle x_1, x_2, x_3, x_4, x_5 ~\mid  ~x_1=x_5x_2^2,
~x_4=x_5^2, ~x_4x_2^3, ~x_1x_5^{-1}x_3 \rangle\\
\cong \langle x_1, x_2, x_3, x_4, x_5, x_6 ~\mid  ~x_1^{-1}x_5x_2^2,
~x_4^{-1}x_5^2, ~x_6^{-1}x_4x_2, ~x_6x_2^2, ~x_1x_5^{-1}x_3 \rangle
= H$.

\end{example}

We note that this procedure of breaking relators into pieces of
length 3 increases the total relator length (measured as the sum of
the length of all relators) by at most a factor of 2.

Since we need our  ``elementary'' isomorphisms to be also given in
the form $x_i \to y_i$, we note that the isomorphism between the
first two presentations above is given by $x_i \to x_i, ~i=1, 2, 3,$
and the inverse isomorphism is given by $x_i \to x_i, ~i=1, 2, 3;
~x_4 \to x_1^2$. By composing elementary isomorphisms, we compute
the isomorphism $\varphi$ between the first and  the last
presentations: $\varphi: x_1 \to x_5, ~x_2 \to x_2, ~x_3 \to x_3$.
By composing the inverses of  elementary isomorphisms, we compute
$\varphi^{-1}: x_1 \to x_1x_2^2, ~x_2 \to x_2, ~x_3 \to x_3, ~x_4
\to x_1^2, ~x_5 \to x_1, ~x_6 \to x_1^2x_2$.  We see that even in
this toy example, recovering $\varphi^{-1}$ from the  public
$\varphi$ is not quite trivial without knowing a sequence of
intermediate Tietze transformations. Furthermore, if Alice discards,
say, two of the relators from the last presentation to get a public
$\hat H = \langle x_1, x_2, x_3, x_4, x_5, x_6 ~\mid
~x_1^{-1}x_5x_2^2,  ~x_6x_2^2, ~x_1x_5^{-1}x_3$, then there is no
isomorphism between $\hat H$ and $G$ whatsoever, and the problem for
the adversary is now even less trivial: to find relators completing
the public presentation $\hat H$ to a presentation  $H$ isomorphic
to $G$ by way of the public isomorphism $\varphi$, and then find
$\varphi^{-1}$.

Moreover, $\varphi$ as a map on the generators of $G$ may not induce
an {\it onto} homomorphism from $G$ to $\hat H$, and this will
deprive the adversary even from the ``brute force" attack by looking
for a map $\psi$  on the generators of $\hat H$ such that  $\psi:
\hat H \to G$ is a homomorphism, and $\psi(\varphi)$ is identical on
$G$. If, say, in the example above we discard the relator
$x_1x_5^{-1}x_3$ from the final presentation $H$, then $x_4$ will
not be in the subgroup of $\hat H$ generated by $\varphi(x_i)$, and
therefore there cannot possibly be a  $\psi: \hat H \to G$ such that
$\psi(\varphi)$ is identical on $G$.

We now describe a homomorphic public key encryption scheme a little
more formally.
\begin{description}
\item[] \textsf{Key Generation:}
  Let $\varphi : G \to H$ be an isomorphism.  Alice's public key then
  consists of $\varphi$ as well as presentations $G$ and $\hat H$,
  where $\hat H$ is obtained from $H$ by keeping all of the generators
  but discarding some of the relators.  Alice's private key consists
  of $\varphi^{-1}$ and $H$.

\item[] \textsf{Encrypt:}
  Bob's plaintext is $g \in G$. To encrypt, he selects an arbitrary word $w_g$
in the generators of $G$ representing the element $g$ and applies
the public isomorphism $\varphi$ to $w_g$ to get $\varphi(w_g)$,
which is a word in the generators of $H$ (or $\hat H$, since $H$ and
$\hat H$ have the same set of generators). He then selects an
arbitrary word $h_g$ in the generators of $\hat H$ representing the
same element of $\hat H$ as $\varphi(w_g)$ does, and this is now his
ciphertext:  $h_g = E(g)$. 

\item[] \textsf{Decrypt:}
  To decrypt, Alice applies her private map
$\varphi^{-1}$  to $h_g$ to get a word $w'_g = \varphi^{-1}(h_g)$.
This word $w'_g$ represents the same element of $G$ as $w_g$ does
because $\varphi^{-1}(h_g) = \varphi^{-1}(\varphi(w_g)) = w_g$ in
the group $G$ since both $\varphi$  and $\varphi^{-1}$ are
homomorphisms, and the composition of $\varphi$  and $\varphi^{-1}$
is the identity map on the group $G$.
\end{description}

\noindent In the following example, we use the
presentations

\noindent  $G=\big\langle x_1, x_2, x_3 ~\mid ~x_1^2x_2^3, ~x_1x_2^2x_1^{-1}x_3
\big\rangle,\quad \hat H= \big\langle x_1, x_2, x_3, x_4, x_5, x_6
~\mid ~x_1=x_5x_2^2, ~x_4=x_5^2,\\
x_6=x_4x_2, ~x_6x_2^2
\big\rangle,$ and the isomorphism $\varphi: x_1 \to
x_5, ~x_2 \to x_2, ~x_3 \to x_3$ from Example \ref{Tietze_example}
to illustrate how encryption works.

\begin{example} \label{encryption_example}
Let the plaintext be the element $g \in G$ represented by the word
$x_1x_2$. Then $\varphi(x_1x_2) = x_5x_2$. Then the word $x_5x_2$ is
randomized in $\hat H$ by using relators of $\hat H$ as well as
``trivial'' relators $x_ix_i^{-1}=1$ and $x_i^{-1}x_i=1$. For
example:  multiply $x_5x_2$ by $x_4x_4^{-1}$ to get
$x_4x_4^{-1}x_5x_2$. Then replace $x_4$ by $x_5^2$, according to one
of the relators of $\hat H$, and get $x_5^2x_4^{-1}x_5x_2$. Now
insert $x_6x_6^{-1}$ between $x_5$ and  $x_2$ to get
$x_5^2x_4^{-1}x_5x_6x_6^{-1}x_2$, and then replace $x_6$ by $x_4x_2$
to get $x_5^2x_4^{-1}x_5x_4x_2x_6^{-1}x_2$, which can be used as the
encryption $E(g)$.

\end{example}

Finally, we note that {\it automorphisms}, instead of general
isomorphisms, were used in \cite{Tropical} and \cite{Moh} to build
public key cryptographic primitives employing the same general idea
of building an automorphism as a composition of elementary ones. In
\cite{Moh}, those were automorphisms of a polynomial algebra, while
in \cite{Tropical} automorphisms of a tropical algebra were used
along the same lines. We also note that ``elementary isomorphisms"
(i.e., Tietze transformations) are universal in nature and can be
adapted to most any algebraic structure, see e.g. \cite{MSYubook},
\cite{Umirbaev}, and \cite{ShYu}.

\section{Semidirect product of groups and  more peculiar computational assumptions}
\label{semidirect}

\noindent Using a semidirect product of (semi)groups as the platform
for a very simple key exchange protocol (inspired by the
Diffie-Hellman  protocol) yields new and sometimes rather peculiar
computational assumptions. The exposition in this section follows
\cite{semidirect1} (see also \cite{semidirect2}).

First we recall the definition of a semidirect product:

\begin{definition} Let $G, H$ be two groups, let $Aut(G)$ be the group of automorphisms of $G$,
and let $\rho: H \rightarrow Aut(G)$ be a homomorphism. Then the
semidirect product of $G$ and $H$ is the set
$$\Gamma = G \rtimes_{\rho} H = \left \{ (g, h): g \in G, ~h \in H \right \}$$
with the group operation given by\\
\centerline{$(g, h)(g', h')=(g^{\rho(h')} \cdot  g', ~h \cdot h')$.}\\
Here $g^{\rho(h')}$ denotes the image of  $g$ under the automorphism
$\rho(h')$, and when we write a product $h \cdot h'$ of two
morphisms, this means that $h$ is applied first.
\end{definition}

In this section, we focus on a special case of this construction,
where the group $H$ is just a subgroup of the group $Aut(G)$. If
$H=Aut(G)$, then the corresponding semidirect product is called the
{\it holomorph} of the group $G$. Thus,  the holomorph of $G$,
usually denoted by $Hol(G)$, is the set of all pairs $(g, ~\phi)$,
where $g \in G, ~\phi \in Aut(G)$, with the group operation given by
~$(g,~\phi)\cdot  (g',~\phi') = (\phi'(g)\cdot g',~\phi \cdot
\phi')$.

It is often more practical to use a subgroup of $Aut(G)$ in this
construction, and this is exactly what we do below, where we describe a key exchange protocol
that uses (as the platform) an extension of a group $G$ by a cyclic
group of automorphisms.

One can also use this construction if $G$ is not necessarily a
group, but just a semigroup, and/or consider endomorphisms of $G$,
not necessarily automorphisms. Then the result will be a semigroup

Thus, let $G$ be a (semi)group. An element $g\in G$ is chosen and
made public as well as an arbitrary automorphism $\phi\in Aut(G)$
(or an arbitrary endomorphism $\phi\in End(G)$). Bob chooses a
private $n\in \mathbb{N}$, while Alice chooses a private $m\in
\mathbb{N}$. Both Alice and Bob are going to work with elements of
the form $(g, \phi^r)$, where $g\in G, ~r\in \mathbb{N}$. Note that
two elements of this form are multiplied as follows:  ~$(g, \phi^r)
\cdot (h, \phi^s) = (\phi^s(g)\cdot h, ~\phi^{r+s})$.
\medskip

The following is a public key exchange protocol between Alice and
Bob.

\begin{enumerate}

\item Alice computes $(g, \phi)^m = (\phi^{m-1}(g) \cdots \phi^{2}(g) \cdot \phi(g) \cdot g,
~\phi^m)$ and sends {\bf only the first component} of this pair to
Bob. Thus, she sends to Bob {\bf only} the element $a =
\phi^{m-1}(g) \cdots \phi^{2}(g) \cdot \phi(g) \cdot g$ of the
(semi)group $G$.
\medskip

\item Bob computes $(g, \phi)^n = (\phi^{n-1}(g) \cdots \phi^{2}(g) \cdot \phi(g) \cdot g,
~\phi^n)$ and sends {\bf only the first component} of this pair to
Alice. Thus, he sends to Alice {\bf only} the element $b =
\phi^{n-1}(g) \cdots \phi^{2}(g) \cdot \phi(g) \cdot g$ of the
(semi)group $G$.
\medskip

\item Alice computes $(b, x) \cdot (a, ~\phi^m) = (\phi^m(b) \cdot a,
~x  \cdot \phi^{m})$. Her key is now $K_A = \phi^m(b) \cdot a$. Note
that she does not actually ``compute" $x \cdot \phi^{m}$ because she
does not know the automorphism $x=\phi^{n}$; recall that it was not
transmitted to her. But she does not need it to compute $K_A$.
\medskip

\item Bob computes $(a, y) \cdot (b, ~\phi^n) = (\phi^n(a) \cdot b,
~y \cdot \phi^{n})$. His key is now $K_B = \phi^n(a) \cdot b$.
Again, Bob does not actually ``compute" $y \cdot \phi^{n}$ because
he does not know the automorphism $y=\phi^{m}$.
\medskip

\item Since $(b, x) \cdot (a, ~\phi^m) = (a, ~y) \cdot (b, ~\phi^n) =
(g, ~\phi)^{m+n}$, we should have $K_A = K_B = K$, the shared secret
key.

\end{enumerate}

\begin{remark}
Note that, in contrast with the original Diffie-Hellman key
exchange, correctness here is based on the equality $h^{m}\cdot
h^{n} = h^{n} \cdot h^{m} =  h^{m+n}$  rather  than on the equality
$(h^{m})^{n} = (h^{n})^{m} = h^{mn}$. In  the original
Diffie-Hellman set up, our trick would not work because, if the
shared key $K$ was just the product of two  openly transmitted
elements, then anybody, including the eavesdropper, could compute
$K$.
\end{remark}

We note that the general protocol above can be used  with {\it any} non-commutative
group $G$  if $\phi$ is selected to be a non-trivial inner
automorphism, i.e., conjugation by an element which is not in the
center of $G$. Furthermore, it can be used  with  any
non-commutative {\it semigroup} $G$ as well, as long as  $G$ has
some invertible elements; these can be used to produce  inner
automorphisms. A typical example of such a semigroup would be a
semigroup of matrices over some ring.

Now let $G$ be a non-commutative (semi)group and let  $h \in G$ be
an invertible non-central element. Then conjugation by $h$ is a
non-identical inner automorphism of $G$ that we denote by
$\varphi_{h}$. We use an extension of the semigroup $G$ by the inner
automorphism $\varphi_{h}$, as described in the beginning of this
section.  For any element $g \in G$ and for any integer $k \ge 1$,
we have

$$\varphi_{h}(g) = g^{-1} g h; ~\varphi^k_{h}(g) = h^{-k} g h^k.$$

\noindent Now our general protocol is specialized in this case as follows.

\medskip

\begin{enumerate}

\item Alice and Bob agree on a (semi)group $G$ and on public elements $g, h \in G$, where $h$ is
an invertible non-central element.
\medskip

\item  Alice selects a private positive integer
$m$, and Bob selects a private positive integer $n$.
\medskip

\item Alice computes $(g, \varphi_{h})^m = (h^{-m+1} g h^{m-1} \cdots h^{-2} g h^2 \cdot h^{-1} g h \cdot g,
~\varphi_{h}^m)$ and sends {\bf only the first component} of this
pair to Bob. Thus, she sends to Bob {\bf only} the element\\
$$A = h^{-m+1} g h^{m-1} \cdots h^{-2} g h^2 \cdot h^{-1} g h \cdot
g = h^{-m} (hg)^{m}.$$
\medskip

\item Bob computes $(g, \varphi_{h})^n = (h^{-n+1} g h^{n-1} \cdots h^{-2} g h^2 \cdot h^{-1} g h \cdot g,
~\varphi_{h}^n)$ and sends {\bf only the first component} of this
pair to Alice. Thus, he sends to Alice {\bf only} the element\\
$$B = h^{-n+1} g h^{n-1} \cdots h^{-2} g h^2 \cdot h^{-1} g h \cdot
g = h^{-n} (hg)^{n}.$$


\item Alice computes $(B, x) \cdot (A, ~\varphi_{h}^m) = (\varphi_{h}^m(B) \cdot A,
~x  \cdot \varphi_{h}^{m})$. Her key is now $K_{Alice} =
\varphi_{h}^m(B) \cdot A = h^{-(m+n)}(hg)^{m+n}$. Note that she does
not actually ``compute" $x \cdot \varphi_{h}^{m}$ because she does
not know the automorphism $x=\varphi_{h}^{n}$; recall that it was
not transmitted to her. But she does not need it to compute
$K_{Alice}$.
\medskip

\item Bob computes $(A, y) \cdot (B, ~\varphi_{h}^n) = (\varphi_{h}^n(A) \cdot B,
~y \cdot \varphi_{h}^{n})$. His key is now $K_{Bob} =
\varphi_{h}^n(A) \cdot B$. Again, Bob does not actually ``compute"
$y \cdot \varphi_{h}^{n}$ because he does not know the automorphism
$y=\varphi_{h}^{m}$.
\medskip

\item Since $(B, x) \cdot (A, ~\varphi_{h}^m) = (A, ~y) \cdot (B, ~\varphi_{h}^n) =
(M, ~\varphi_{h})^{m+n}$, we should have $K_{Alice} = K_{Bob} = K$,
the shared secret key.

\end{enumerate}

Thus, the shared secret key in this protocol is
$$K = \varphi_{h}^m(B) \cdot A = \varphi_{h}^n(A) \cdot B = h^{-(m+n)}(hg)^{m+n}.$$

Therefore, our security assumption here is that it is
computationally hard to retrieve the key $K = h^{-(m+n)}(hg)^{m+n}$
from the quadruple
$(h, ~g, ~h^{-m}(hg)^{m}, ~h^{-n}(hg)^{n})$. In particular, we have
to take care that the elements $h$ and  $hg$ do not commute because
otherwise, $K$ is just a product of $h^{-m}(hg)^{m}$  and
$h^{-n}(hg)^{n}$. Once again, the problem is:

\begin{quote}
Given a (semi)group $G$ and elements $g, ~h, ~h^{-m}(hg)^{m}$, and
$h^{-n}(hg)^{n}$ of $G$, find $h^{-(m+n)}(hg)^{m+n}$.

\end{quote}

Compare this to the Diffie-Hellman problem from Section \ref{DH}:

\begin{quote}
 Given a (semi)group $G$ and elements $g, g^{n}$, and  $g^{m}$ of
 $G$, find $g^{mn}$.

\end{quote}

A weaker security assumption arises if an eavesdropper tries to
recover a private exponent from a transmission, i.e., to recover,
say, $m$ from $h^{-m} (hg)^{m}.$ A special case of this problem,
where $h=1$, is the ``discrete log" problem, namely: recover $m$
from $g$ and $g^{m}$. However, the ``discrete log" problem is a
problem on cyclic, in particular abelian, groups, whereas in the
former problem it is essential that $g$ and $h$ do not commute.

By varying the automorphism (or endomorphism) used for an extension
of $G$, one can get many other security assumptions. However, many
(semi)groups $G$ just do not have outer (i.e., non-inner)
automorphisms, so there is no guarantee that a selected platform
(semi)group will have any outer automorphisms. On the other hand, it
will have inner automorphisms as long as it has invertible
non-central elements.

In conclusion, we note that there is always a concern (as well as in
the standard Diffie-Hellman protocol) about the orders of  public
elements (in our case, about the orders of $h$ and $hg$): if one of
the orders is too small, then a brute force attack may be feasible.

If a group of matrices of small size is chosen as the platform, then
the above protocol turns out to be vulnerable to a ``linear algebra
attack", similar to an attack on Stickel's protocol
\cite{Stickel_self} offered in \cite{Stickel}, albeit more
sophisticated, see \cite{MR}, \cite{R}, \cite{Romankov:2016}. A
composition of conjugating automorphism with a field automorphism
was employed in \cite{KLS}, but this automorphism still turned out
to be not complex enough to make the protocol withstand a linear
algebra attack, see \cite{DMU}, \cite{R}. Selecting a good platform
(semi)group for the protocol in this section still remains an open
problem.

Finally, we mention another, rather different, proposal \cite{PHKCP}
of a cryptosystem based on the semidirect product of two groups and
yet another, more complex, proposal of a key agreement based on the
semidirect product of two monoids \cite{AAGL}.

\section{The subset sum and the knapsack problems}
\label{knapsack}

As usual, elements of a group $G$ are given as words in the alphabet
$X \cup X^{-1}$. We begin with three decision problems:

\medskip
\noindent{\bf The subset sum  problem ($\SSP$):} Given
$g_1,\ldots,g_k,g\in G$ decide if
  \begin{equation} \label{eq:SSP-def}
  g = g_1^{\varepsilon_1} \ldots g_k^{\varepsilon_k}
  \end{equation}
for some $\varepsilon_1,\ldots,\varepsilon_k \in \{0,1\}$.

\medskip
\noindent {\bf The knapsack problem ($\KP$):}  Given
$g_1,\ldots,g_k,g\in G$ decide if
\begin{equation}\label{eq:IKP-def}
g = g_1^{\varepsilon_1} \ldots g_k^{\varepsilon_k}
\end{equation}
for some  non-negative integers
$\varepsilon_1,\ldots,\varepsilon_k$.

The third problem is equivalent to $\KP$ in the abelian case, but in
general this is a completely different problem:

\medskip
\noindent {\bf The  Submonoid membership problem ($\mathbf{SMP}$)}:
Given elements $g_1,\ldots,g_k,g\in G$ decide if $g$ belongs to the
submonoid generated by $g_1, \ldots, g_k$ in $G$, i.e.,  if the
following equality holds for some $g_{i_1}, \ldots, g_{i_s} \in
\{g_1, \ldots, g_k\}, s \in \mathbb{N}$:
\begin{equation}\label{eq:SMP-def}
g = g_{i_1}, \ldots, g_{i_s}.
\end{equation}

\medskip
The restriction of $\SMP$ to the case where the set of generators
$\{g_1, \ldots,g_n\}$ is closed under inversion (so that the
submonoid is actually a subgroup of $G$) is  a well-known {\em
subgroup membership problem}, one of the most basic algorithmic
problems in group theory.

\medskip
There are also natural {\it search} versions of the decision
problems above, where the goal is to find a particular solution to
the equations (\ref{eq:SSP-def}), (\ref{eq:IKP-def}), or
(\ref{eq:SMP-def}), provided that solutions do  exist.

We also mention, in passing, an interesting research avenue explored
in \cite{knapsack}: many search problems can be converted to {\it
optimization} problems asking for an ``optimal" (usually meaning
``minimal") solution of the corresponding search problem. A
well-known example of an optimization problems is the {\it geodesic
problem}: given a word in the generators of a group $G$, find a word
of minimum length representing the same element of $G$.

The classical (i.e., not group-theoretical) subset sum problem is
one of the very basic $\NP$-complete problems, so there is extensive
related bibliography (see \cite{Kellerer-Pferschy-Pisinger:2004}).
The $\SSP$ problem attracted a lot of extra attention when Merkle
and Hellmann designed a public key cryptosystem \cite{Merkle} based
on a variation of $\SSP$. That cryptosystem was broken by Shamir in
\cite{Shamir:1984}, but the interest persists and the ideas survive
in numerous  new cryptosystems and their variations (see e.g.
\cite{Odlyzko:1990}). Generalizations of knapsack-type cryptosystems
to non-commutative groups seem quite promising from the viewpoint of
post-quantum cryptography, although relevant cryptographic schemes
are yet to be built.

In \cite{knapsack}, the authors showed, in particular, that $\SSP$
is $\NP$-complete in: (1) the direct sum of countably many  copies
of the infinite cyclic group $\Z$; (2) free metabelian  non-abelian
groups of finite rank; (3) wreath product of two finitely generated
infinite abelian groups; (4) Thompson's group $F$; (5)
Baumslag--Solitar group $BS(m,n)$ for $|m|\ne |n|$, and in many
other groups.

In \cite{KLZ}, the authors showed that the subset sum problem is
polynomial time decidable in every finitely generated virtually
nilpotent group but there exists a polycyclic group where this
problem is NP-complete. Later in  \cite{polycyclic Subset sum},
Nikolaev and Ushakov showed that, in fact, every polycyclic
non-virtually-nilpotent group has NP-complete subset sum problem.

Also in \cite{KLZ}, it was shown that the knapsack problem is
undecidable in a direct product of sufficiently many copies of the
discrete Heisenberg group (which is nilpotent of class 2). However,
for the discrete Heisenberg group itself, the knapsack problem is
decidable. Thus,  decidability of the knapsack problem is not
preserved under direct products. In \cite{Frenkel}, the effect of
free and direct products on the time complexity of the knapsack and
related problems was studied further.

\section{The Post correspondence problem}
\label{Post}

The Post correspondence problem  $\PCP(\CA)$ for a semigroup (or any
other algebraic structure) $\CA$  is to decide, given two n-tuples
$u = (u_1, \ldots,u_n)$ and $v = (v_1, \ldots,v_n)$ of elements of
$\CA$,  if there is a term (called a {\em solution}) $t(x_1,
\ldots,x_n)$  in the language of $\CA$ such that $t(u_1, \ldots,u_n)
= t(v_1, \ldots,v_n)$ in $\CA$. In 1946 Post introduced this problem
in the case of free monoids (free semigroups)  and proved that it is
undecidable \cite{Post:1946}. (See \cite{Sipser:2005} for a simpler
proof.)

The $\PCP$ in groups is closely related to the problem  of finding
the {\em  equalizer} $E(\phi,\psi)$ of two group homomorphisms
$\phi, \psi : H \to G$.   The equalizer  is defined as $E(\phi,\psi)
= \{w \in H \mid \phi(w) = \psi(w)\}$.  Specifically, $\PCP$ in a
group $G$ is the same as to decide if the equalizer of a given pair
of homomorphisms $\phi, \psi \in Hom(H,G)$, where $H$ is a free
group of finite rank in the variety $\Var(G)$ generated by $G$,  is
trivial or not. Indeed, in this case every tuple $u = (u_1, \ldots,
u_n)$ of elements of $G$ corresponds to a homomorphism $\phi_u$ from
a free group $H$ with a basis $x_1, \ldots, x_n$ in the variety
$\Var(G)$ such that $\phi_u(x_1) =  u_1, \ldots, \phi_u(x_n) = u_n$.
The equalizer $E(\phi_u,\phi_v)$ describes all solutions $w$ for the
instance $(u, v)$.

There is an interesting variation of the Post correspondence problem
for semigroups and groups that we call a non-homogeneous Post
correspondence problem, or a {\it general} Post correspondence
problem $\GPCP$, following \cite{Post}: given two tuples $u$ and $v$
of elements in a (semi)group $S$ as above and two extra elements $a,
b \in S$, decide if there is a term $t(x_1, \ldots,x_n)$ such that
$a t(u_1, \ldots,u_n) = b t(v_1, \ldots,v_n)$ in $S$. Interesting
connections between $\GPCP$ and the (double) twisted conjugacy
problem were reported in \cite{Post}.  Specifically, it was shown in
\cite{Post} that the double endomorphism twisted conjugacy problem
in a relatively free group in $\Var(G)$ is equivalent to $\GPCP(G)$,
and, in general, the double endomorphism twisted conjugacy problem
in $G$ is $\P$-time reducible to $\GPCP(G)$.

Another interesting observation made in \cite{Post} is that if
$\GPCP$ is decidable in a group $G$ then there is a {\em uniform}
algorithm to solve the word problem in every finitely presented
(relative to $G$) quotient of $G$. Furthermore, since decidability
of $\GPCP$ in $G$ is inherited by all subgroups of $G$, decidability
of $\GPCP$ in $G$ implies the uniform decidability of the word
problem in every finitely presented quotient of every subgroup of
$G$.

Examples of groups with undecidable $\GPCP$ include free groups  and
free solvable groups of derived length at least 3 and  sufficiently
high rank (\cite{Post}). On the other hand, examples of groups where
$\GPCP$ is decidable in polynomial time include all finitely
generated nilpotent groups.

Furthermore,
it was shown in \cite{Post} that in a free group,  the bounded
$\GPCP$ is $\NP$-complete. (In the bounded version of $\GPCP$ one is
looking only for solutions (i.e., the words $t(x_1, \ldots,x_n)$), whose
length is bounded by a given number.)

The {\em search} version of the Post correspondence problem (or of
the bounded version thereof) is to find a solution for a given
instance, provided at least one solution exists. As usual, search
versions can potentially be used to build cryptographic primitives,
although it is not immediately clear how to convert the search
version of the (bounded or not) Post correspondence problem   to a
one-way function with trapdoor.

\section{The hidden subgroup problem}
\label{HSP}

Given a group $G$, a subgroup $H \le G$, and a set $X$, we say that
a function $f : G \to X$ hides the subgroup $H$ if for all $g_1, g_2
\in G$, one has $f(g_1) = f(g_2)$ if and only if $g_1H = g_2 H$ for
the cosets of $H$. Equivalently, the function $f$ is constant on the
cosets of $H$, while it is different between the different cosets of
$H$.

The {\it hidden subgroup problem} (HSP) is:

\begin{quote}
Let $G$ be a finite group, $X$ a finite set, and $f : G \to X$  a
function that hides a subgroup $H \le G$. The function $f$ is given
via an oracle, which uses $O(\log |G|+\log|X|)$ bits. Using
information gained from evaluations of $f$ via its oracle, determine
a generating set for $H$.
\end{quote}

\noindent A special case is where $X$ is a group and $f$ is a group
homomorphism, in which case $H$ corresponds to the kernel of $f$.

The importance of the hidden subgroup problem is due to the facts that:

\begin{itemize}
\item
Shor's polynomial time quantum algorithm for factoring and discrete
logarithm problem (as well as several of its extensions) relies on
the ability of quantum computers to solve the HSP for finite abelian
groups. Both factoring and discrete logarithm problem are of
paramount importance for modern commercial cryptography.

\item  The existence of efficient quantum algorithms for HSP for certain
non-abelian groups would imply efficient quantum algorithms for two major problems:
the graph isomorphism problem and certain shortest vector problems in lattices.
More specifically, an efficient quantum algorithm for the HSP for the symmetric
group would give a quantum algorithm for the graph isomorphism, whereas an efficient
quantum algorithm for the HSP for the dihedral group would give a quantum
algorithm for the shortest vector problem.

\end{itemize}

We refer to \cite{WangWang} for a brief discussion on how the HSP can be generalized to infinite groups.

\section{Relations between some of the problems}
\label{relations}

In this section, we discuss relations between some of the  problems
described earlier in this survey. In the preceding Sections
\ref{knapsack} and \ref{Post} we have already pointed out some of
the relations, now here are some other relations, through the prizm
of cryptographic applications.

We start with the {\it conjugacy search  problem} (CSP), which was
the subject of Section \ref{conjugacy}, and one of its ramifications, the {\it
subgroup-restricted} conjugacy search  problem:

\begin{quote}
Given two elements $w, h$ of a group $G$, a subgroup  $A \leq G$,
and the information that $w^a=h$ for some $a \in A$, find at least
one particular element $a$ like that.
\end{quote}

In reference to the Ko-Lee protocol described in Section \ref{conjugacy},
one of the parties (Alice) transmits $w^a$ for some private $a \in
A$, and the other party (Bob) transmits $w^b$ for some private $b
\in B$, where the subgroups $A$ and  $B$ commute elementwise, i.e.,
$ab=ba$ for any  $a \in A$,  $b \in B$.

Now suppose the adversary finds $a_1, a_2 \in A$ such that
$a_1 w a_2 = a^{-1} w a$  and  $b_1, b_2 \in B$ such that $b_1 w b_2
= b^{-1} w b$.  Then the adversary  gets
$$a_1 b_1 w b_2 a_2 =
a_1 b^{-1} w b a_2 = b^{-1} a_1 w a_2 b = b^{-1} a^{-1} w a b =K,$$

\noindent the shared secret key.

We emphasize that these $a_1, a_2$  and $b_1, b_2$ do not have
anything to do with the private elements originally selected by
Alice or Bob, which simplifies the search substantially. We also
point out that, in fact, it is sufficient for the adversary to find
just one pair, say, $a_1, a_2 \in A$, to get the shared secret  key:
$$a_1 (b^{-1} w b) a_2 = b^{-1} a_1 w a_2 b =b^{-1} a^{-1} w a b =K.$$

 In summary, to get the secret key $K$, the adversary does not have to solve
the (subgroup-restricted) conjugacy search  problem, but instead, it
is sufficient to solve an apparently easier  (subgroup-restricted)
{\it decomposition search  problem}, see  our Section \ref{decomposition}.



Then, one more trick reduces the  decomposition search problem to a
special case where $w=1$, i.e., to the {\it factorization problem},
see our Section \ref{factorization}. Namely, given $w' = a\cdot w
\cdot b$, multiply it on the left by the  element $w^{-1}$ (which is
the inverse of the public element $w$) to get
    $$w'' = w^{-1}a\cdot w\cdot b = (w^{-1}a\cdot w)\cdot b.$$

Thus, if we denote by $A^w$ the subgroup conjugate to $A$ by the (public)
element $w$, the problem for the adversary is now the following {\it factorization search
problem}:

\begin{quote}
 Given  an element $w'$ of a group $G$ and two subgroups
$A^w, B \leq G$, find any two elements $a \in A^w$ and $b \in B$
that would   satisfy $a\cdot  b = w'$, provided at least one such
pair of elements exists.
\end{quote}

 Since in the original Ko-Lee protocol one has $A=B$, this yields
 the following interesting observation: if in that protocol $A$ is
a normal subgroup of $G$, then $A^w=A$, and the above problem
becomes: given $w' \in A$, find any two elements $a_1, a_2  \in  A$
such that $w' = a_1 a_2$. This problem is trivial: $a_1$ here could
be any element from $A$, and then $a_2 = a_1^{-1} w'$.

 Therefore, in choosing the platform group $G$ and two commuting
subgroups for a protocol described in our Section \ref{conjugacy} or
Section \ref{decomposition}, one has to avoid normal subgroups. This
means, in particular, that ``artificially" introducing commuting
subgroups as, say, direct factors is inappropriate from the security
point of view.

At the other extreme, there are  {\it malnormal} subgroups. A
subgroup $A \leq G$ is called malnormal in $G$ if, for any $g \in
G$, $A^g \cap A = \{1\}$. We observe that if, in the original Ko-Lee
protocol, $A$ is a malnormal subgroup of $G$, then the decomposition
search problem corresponding to that protocol has a unique solution
if $w \notin A$. Indeed, suppose $w' = a_1\cdot w \cdot a'_1 =
a_2\cdot w \cdot a'_2$, where $a_1 \ne a_2$, say. Then $a_2^{-1} a_1
w = w a'_2 a_1'^{-1}$, hence $w^{-1} a_2^{-1} a_1 w = a'_2
a_1'^{-1}$. Since $A$ is  malnormal, the element on the left does
not belong to $A$, whereas the one on the right does, a
contradiction. This argument shows that, in fact, already if $A^w
\cap A = \{1\}$ for this particular $w$, then the corresponding
decomposition search problem has a unique solution.

Finally, we describe one more trick that reduces, to some extent,
the decomposition search  problem to the (subgroup-restricted)
conjugacy search problem.  Suppose we are given $w' = aw b$, and we need to recover
$a \in A$ and $b \in B$, where $A$ and $B$ are two elementwise
commuting subgroups of a group $G$.

 Pick any $b_1 \in B$ and compute:
$$[awb, b_1] = b^{-1}w^{-1}a^{-1}b_1^{-1}awbb_1 = b^{-1}w^{-1}b_1^{-1}wbb_1
= (b_1^{-1})^{wb} b_1 = ((b_1^{-1})^{w})^b b_1.$$

Since we know $b_1$, we can multiply the result by $b_1^{-1}$ on the
right to get  $w'' = ((b_1^{-1})^{w})^b$. Now the problem becomes:
recover $b \in B$ from the known  $w'' = ((b_1^{-1})^{w})^b$ and
$(b_1^{-1})^{w}$. This is the subgroup-restricted conjugacy search
problem. By solving it, one can recover a $b \in B$.

 Similarly, to recover an $a \in A$, one picks any $a_1 \in A$ and computes:
\begin{align*}
[(awb)^{-1}, ~(a_1)^{-1}] &= awba_1b^{-1}w^{-1}a^{-1}a_1^{-1}\\
&= awa_1w^{-1}a^{-1}a_1^{-1}=(a_1)^{w^{-1}a^{-1}} a_1^{-1} =
((a_1)^{w^{-1}})^{a^{-1}} a_1^{-1}.
\end{align*}

 Multiply the result by $a_1$ on the right to get  $w'' =
((a_1)^{w^{-1}})^{a^{-1}}$, so that the problem becomes: recover $a
\in A$ from the known  $w'' = ((a_1)^{w^{-1}})^{a^{-1}}$ and
$(a_1)^{w^{-1}}$.

 We have to note  that, since a solution of the subgroup-restricted conjugacy search
problem is not always unique, solving the above two instances of
this problem may not necessarily give the right solution of the
original decomposition problem. However, any two solutions, call
them $b'$  and $b''$, of the first conjugacy search problem  differ
by an element of the centralizer of $(b_1^{-1})^{w}$, and this
centralizer is unlikely to have a non-trivial intersection with $B$.

A similar computation shows that the  same trick   reduces the
factorization search problem, too, to the subgroup-restricted
conjugacy search problem. Suppose we are given $w' = a b$, and we
need to recover $a \in A$ and $b \in B$, where $A$ and $B$ are two
elementwise commuting subgroups of a group $G$.  Pick any $b_1 \in
B$ and compute
$$[ab, b_1] = b^{-1}a^{-1}b_1^{-1}a b b_1 = (b_1^{-1})^{b} b_1.$$

Since we know $b_1$, we can multiply the result by $b_1^{-1}$ on the
right to get  $w'' = (b_1^{-1})^{b}$. This is the
subgroup-restricted conjugacy search problem. By solving it, one can
recover a $b \in B$.

 This same trick can, in fact, be used to attack the subgroup-restricted conjugacy search
problem itself. Suppose we are given $w' = a^{-1} w a$, and we need
to recover $a \in A$. Pick any $b$ from the centralizer of $A$;
typically, there is a public subgroup $B$ that commutes with $A$
elementwise; then just pick any $b \in B$. Then compute
$$[w', b] = [a^{-1} w a, b] = a^{-1} w^{-1} a b^{-1} a^{-1} wab = a^{-1} w^{-1} b^{-1}wab =
(b^{-w})^a b.$$

 Multiply the result by $b^{-1}$ on the right to get  $(b^{-w})^a$,
so the problem now is to recover $a \in A$ from $(b^{-w})^a$ and
$b^{-w}$. This problem might be easier than the original problem
because there is flexibility in choosing $b \in B$. In particular, a
feasible attack might be to  choose several different $b \in B$  and
try to solve the above conjugacy search problem for each in parallel
by using some general method  (e.g., a length-based attack). Chances
are that the attack will be successful for at least one of the
$b$'s.

\end{document}